\documentclass[a4paper,10pt]{article}
\usepackage{amsmath, amsthm, amssymb}
\usepackage{graphicx}

\title{PL-virtual knots}
\author{Neil R. Nicholson\\
	William Jewell College\\
	nicholsonn$@$william.jewell.edu\\
	}

\newtheoremstyle{dotless}{}{}{\itshape}{}{\bfseries}{}{ }{}
\def\Real{\hbox{I\kern-.1667em\hbox{R}}}
\theoremstyle{dotless}	
\input epsf

\newtheorem{thm}{Theorem}[section]

\newtheorem{lem}[thm]{Lemma}

\newcommand{\kb}[1]{\langle #1 \rangle}

\begin{document}

\maketitle

\begin{abstract}
Piecewise-linear virtual knots are discussed and classified up to edge index six.  
\end{abstract}

\section{Introduction}
\label{intro}

The piecewise-linear approach to knots has been utilized and studied for quite some time, and it is still the topic of numerous current research projects and papers.  While the major question in the field is simple to state (What is the minimum number of edges needed to construct a given knot?), relatively little is known about its answer.  

Richard Randell \cite{RandellInvariants} took some of the first steps towards understanding the space of piecewise linear knots and its relationships to certain knot invariants.  Specifically, he found the minimal edge number for all knots with six or fewer crossings.  His ideas for projections of PL-knots have also proven to be quite helpful.  Since his paper, numerous other advances in the field have occurred, including but not limited to \cite{AdamsBrennan}, \cite{Calvo}, \cite{McCabe}, and \cite{Meissen}.

\begin{figure}[hbtp]
  \begin{center}
    \leavevmode
    \epsfxsize = 10cm
    \epsfysize = 1.8cm
    \epsfbox{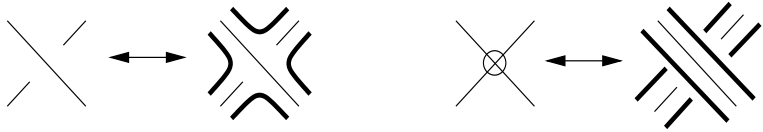}
    \caption{Classical and virtual crossings in surfaces, respectively}
    \label{surfaces}
  \end{center}
\end{figure}

In 1999, Louis Kauffman \cite{KauffmanVirtual} formally introduced to the mathematical community the theory of virtual knots.  While classical knots lie in $\mathbb{R}^{3}$ or $\mathbb{S}^{3}$, virtual knots are embeddings of curves in thickened compact, connected surfaces of positive genus \cite{Kuperberg}.  Figure \ref{surfaces} depicts a real crossing and a virtual crossing in such a space.   Though still a relatively new topic, virtual knot theory has been thoroughly studied and applications to other areas (though the field is deserving enough of its own study) of topology abound.  Interesting questions and directions for research can be found in \cite{KauffmanUnsolved}.

This paper is inspired by the above two ideas and the realization that there is indeed a nontrivial virtual knot with five edges, though the first nontrivial classical knot requires six edges.  Specifically, in Section \ref{defs}, we lay the necessary groundwork and show that Kauffman's ice model for the Jones polynomial \cite{KauffmanStatistical} carries over to virtual knots.  Section \ref{edges} considers the space of six-edged virtual knots and finishes with our main theorem, a classification of all such virtual knots.  We finish with some questions for future research.

\section{Definitions and Preliminaries}
\label{defs}
\begin{figure}[hbtp]
  \begin{center}
    \leavevmode
    \epsfxsize = 10cm
    \epsfysize = 4cm
    \epsfbox{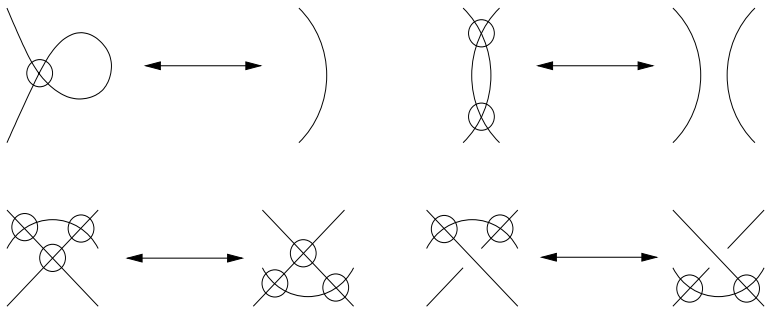}
    \caption{Virtual Reidemeister moves}
    \label{virtualreidemeister}
  \end{center}
\end{figure}

A \textit{virtual knot diagram} differs from a classical knot diagram only by \textit{virtual crossings}.  The standard picture for a virtual crossing is a $4$-valent vertex enclosed in a small circle, as pictured in Fig. \ref{virtualreidemeister}.  We think of these crossings as \textit{not being crossings} and therefore not existing in the real sense (hence the term ``virtual").  A \textit{virtual knot} is an equivalence class of diagrams where equivalence is determined by the classical Reidemeister moves and four generalized Reidemeister moves involving virtual crossings.  See Fig. \ref{virtualreidemeister}.

By \textit{virtual knot} we will mean a knot which is not equivalent to a knot with no virtual crossings.  A knot with such an equivalence will be called a \textit{classical knot} or \textit{real knot}.   On the other hand, a \textit{virtual diagram} is any diagram with at least one virtual crossing, even if the diagram is equivalent to a classical knot.

Let $D$ be a virtual knot diagram.  The \textit{crossing number} $c(D)$ \textit{of $D$} is the total number of crossings, virtual or real, of the diagram.  Let $c(K)$, the \textit{crossing index} of a virtual knot $K$, be the minimum such value over all diagrams of $K$.  Note that this notation is different than that of the virtual knot table \cite{VirtualTable}.  For a given diagram $D$, one can write $c(D) = c_{v}(D) + c_{r}(D)$, where $c_{v}(D)$ (resp. $c_{r}(D)$) is the number of virtual (resp. real) crossings in $D$.   We define $c_{v}(K)$ to be the minimal number of virtual crossings over all diagrams of a knot type $K$.

We are interested in piecewise-linear virtual knots.  The shadow (projection with classical crossing information removed) of an $n$-edged PL-knot (real or virtual) is called an \textit{$n$-universe} and the number of intersections of an $n$-universe is the number of intersections (we can assume a universe is in general position) between nonadjacent edges of the universe.   The intersection of two adjacent edges is a \textit{vertex} of the diagram.

The \textit{edge index} $e(K)$ of a PL-virtual knot $K$ is the least number of edges required to realize $K$.  Randell \cite{RandellInvariants} has given conditions for an ordered $n$-tuple of points $\lbrace v_{1}$, $v_{2}$, ..., $v_{n}\rbrace$ in $\mathbb{R}^{3}$ to \textit{not} form a classical PL-knot, and those requirements carry over to PL-virtual knots.  If any of the following occur, then a PL-virtual knot is not formed:

\begin{enumerate}
\item $v_{i} = v_{j}$ for some $i \neq j$
\item $v_{i}$ lies on the segment between $v_{j}$ and $v_{k}$, where $i \neq j$, $i \neq k$.
\item Two distinct edges partially intersect or intersect at a point other than a vertex
\end{enumerate}

Note that an ordered $n$-tuple of points in $\mathbb{R}^{3}$ can be considered to form a $PL$-virtual knot by stating when a crossing of edges is a virtual crossing.  The study of the configuration space of these $n$-tuples is a very interesting question and will be the subject of future work.

Negami's proof for the bound on the edge number of classical knots \cite{Negami} carries over to virtual knots, giving:

\begin{thm}
For any PL-knot $K$, $e(K) \leq 2c(K)$.
\end{thm}

To help distinguish virtual knots, though it is not nearly a complete invariant (that is, there are nonequivalent knots with the same Jones polynomial), we utilize the Jones polynomial.  Just as with real knots, there is a state sum formulation of the Jones polynomial via the Kauffman bracket polynomial \cite{KauffmanVirtual}, constructed below.  

Any real crossing of a diagram locally separates the plane into regions, as follows.  If we rotate the overcrossing strand counterclockwise, the regions swept out will be referred to as the \textit{A} regions.  Rotating the strand clockwise sweeps out the $B$ regions.  An $A$-smoothing of the crossing results by locally replacing the crossing with two smooth edges so that the $A$ regions are connected.  We also have $B$-smoothings.

A \textit{state} of $D$ is a choice of smoothing for each real crossing.  Let $a(S)$ and $b(S)$ be the number of $A$- and $B$-smoothings in the state $S$, respectively.

For a diagram $D$, the \textit{Kauffman bracket polynomial} is a Laurent polynomial in the variable $A$ with integer coefficients.  It is given by: 

\begin{center}
$\kb{D} = \sum A^{a(S)-b(S)} (-A^{2} - A^{-2})^{|S|-1}$,
\end{center} 

\noindent where the sum is taken over all states $S$ and $|S|$ is the number of simple closed circuits (i.e., closed loops that ``cross" only via virtual crossings) gotten by smoothing each real crossing according to $S$.  The Kauffman bracket is an invariant of regular isotopy (invariant under real and virtual Reidemeister II and III moves).  If $D$ is oriented, assign a value of $+1$ or $-1$ to each real crossing according to the usual right-hand rule.  Define the writhe of $D$ $w(D)$ to be the sum of these values.  Multiplication by a factor of $(-A)^{-3w(D)}$ gives a version of the Jones polynomial $V_{K}(t)$, an ambient isotopy invariant (invariance under all Reidemeister moves).  The original Jones polynomial results via the substitution $A = t^{-1/4}$.

In short, the Kauffman bracket polynomial is computed for virtual knots exactly as it is computed for real knots by leaving all virtual crossings alone and instead of counting simple closed curves in a state we count the number of simple closed circuits.

\subsection{Ice model}
\label{ice}
We adapt Kauffman's ice model \cite{KauffmanStatistical} for the Jones polynomial for piecewise-linear knots to a model for PL virtual knots.  

Let $K$ be a piecewise-linear virtual knot with diagram $D$ and universe $U$.  The vertices of $D$ correspond to $2$-valent vertices of $U$ while real crossings in $D$ correspond to $4$-valent vertices in $U$.  Virtual crossings in $D$ are inherited as virtual crossings of $U$.  An \textit{arrow covering of $U$} is a choice of orientation for each edge of $U$ such that orientation is preserved through both $2$-valent and virtual vertices and at any $4$-valent vertex, exactly two edges are directed into the vertex and exactly two are directed away from the vertex.  

\begin{figure}[hbtp]
  \begin{center}
    \leavevmode
    \epsfxsize = 5cm
    \epsfysize = 2cm
    \epsfbox{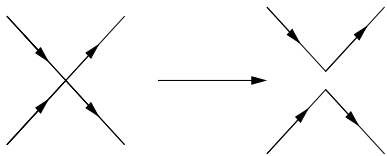}
    \caption{Splitting $4$ valent vertices}
    \label{vertexsplitting}
  \end{center}
\end{figure}

Define a \textit{splitting} of an arrow covered $4$-valent vertex to be a local replacement of the $4$-valent vertex by two $2$-valent vertices such that the arrow covering of the $4$-valent vertex carries over to be an arrow covering of the resultant figure, as in Fig. \ref{vertexsplitting}.  Note that of the six possible local coverings of a $4$-valent vertex, two can split two different ways.  Considering all possible splits, then, an arrow covering corresponds to a collection of diagrams of oriented piecewise-linear closed circuits that may overlap via virtual crossings.

\begin{figure}[hbtp]
  \begin{center}
    \leavevmode
    \epsfxsize = 3cm
    \epsfysize = 2cm
    \put(45,2){$\theta(v)$}
    \epsfbox{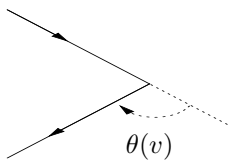}
    \caption{$2$-valent angular contributions}
    \label{2valentcontribution}
  \end{center}
\end{figure}

An oriented $2$-valent vertex $v$ contributes a power of $z$, $z^{\theta(v)}$, via the rule pictured in Fig. \ref{2valentcontribution}.  We make $\theta(v)$ positive for counterclockwise rotation and negative for clockwise rotation.
\vspace{2cm}

If $A$ is an arrow covering of $U$, to each vertex $v$ of $U$ define $K_{A}(v)$ as follows:
\begin{enumerate}
\item If $v$ is $2$-valent, then $K_{A}(v) = z^{\theta(v)}$.
\item If $v$ is $4$-valent, then say $v_{1}$ and $v_{2}$ are the $2$-valent vertices resulting from splitting $v$ via $A$.  Then, let $K_{A}(v)$ be the sum over all possible splittings of $A^{i(v)}z^{\theta(v_{1}) + \theta(v_{2})}$, where $i(v) = 1$ if the split corresponds to an $A$-smoothing of the crossing or $i(v) = -1$ if the split corresponds to a $B$-smoothing of the crossing..
\end{enumerate}

In the case of a real knot, the product over all angular contributions counts the number of oriented simple closed curves in a given state, and Kauffman proved that:
 \begin{center}
 $[K] = \sum \prod K_{A}(v)$,
 \end{center}
where the sum is taken over all arrow coverings $A$ of $U$ and the product is over all vertices $v$.  The same theorem applies in our case, however, the product over angular contributions counts the number of oriented closed circuits of a state.  Thus, we have proven that the ice model carries over to piecewise-linear virtual knots.

\section{Low edge numbers}
\label{edges}
We begin this section by proving a result for PL-virtual knots of any edge number.  Our result mirrors Randell's proposition for real piecewise-linear knots \cite{RandellInvariants}.  

For classical PL-knots, there is a projection in which a given edge $e_{1}$ is not involved in any crossings.  Moreover, the two neighboring edges of $e_{1}$ do not cross.  Randell terms this a \textit{convenient projection}.   Virtual knots, however, do not have such a projection, as the given edges may be involved in virtual crossings.  For calculatory purposes, then, we cannot assume such results about any of the edges of a PL-virtual knot.

\begin{thm}
If $K$ is a PL-virtual knot with $n$ edges, then:
\begin{enumerate}
\item  $c(K) \leq n(n-3)/2$, if $n$ is odd
\item   $c(K) \leq n(n-3)/2 - 1$, if $n$ is even.
\end{enumerate}
\end{thm} \label{edgebound}
\begin{proof}
If $D$ is a diagram of $K$, then any edge cannot cross itself or either of the two edges adjacent to it.  If $n$ is even, it is not possible for each edge to cross the $n-3$ other edges.
\end{proof}

Randell's bound, however, does carry over for certain types of virtual knots:

\begin{thm}
Randell's bound applies to any $n$-edged PL-knot $K$ with \linebreak $n \geq 6c_{v}(K) + 1$.  That is, for any such $K$,
\begin{enumerate}
\item $c(K) \leq (n-1)(n-4)/2$, if $n$ is even. 
\item $c(K) \leq (n-1)(n-4)/2 - (n-3)/3$, if $n$ is odd.
\end{enumerate}
\end{thm}
\begin{proof}
Let $K$ be such a knot.  Then in any projection of $K$ there are at least three consecutive edges not involved in any virtual crossings.  We can then apply Randell's definition of a \textit{convenient projection}: project perpendicular to the middle edge and tilt slightly so that this edge is not involved in any crossings and the first and third of these edges do not cross.
\end{proof}

\subsection{$e(K) \leq 5$}
\label{lessthan5}

A virtual knot diagram $D$ has two \textit{mirror images}.  The \textit{horizontal} mirror image is obtained from $D$ by reflecting the virtual diagram in a mirror and the \textit{vertical} mirror image to $D$ is obtained by switching the over- and undercrossing strands in all real crossings.

In this section we will prove the following:

\begin{thm}
Up to mirror images, the only PL-virtual knots realizable with five or fewer edges are the unknot, $2.1$, $3.5$, and $3.7$. \label{fiveedgesorless}
\end{thm}
\begin{proof}
Assume that $K$ is a PL-virtual knot.  If $e(K) \leq 4$, then by Thm. \ref{edgebound}, $c(K) \leq 1$.  Such knots are unknotted.  Hence we can assume $e(K) = 5$ and $c(K) \leq 5$.  Since $K$ is virtual, we also have that $c_{v}(K) \geq 1$.  We will consider five edges and construct all possible virtual knots.  The one nontrivial virtual knot of crossing index three, $2.1$, is realizable using five edges, as pictured in Fig. \ref{twoone}.  Notice that it not possible to have a universe with four $4$-valent vertices for a knot with five edges.  Therefore we can assume $c(K) = 5$.  Inspection of the virtual knot table shows $c_{v}(K) \leq 2$, yielding two cases.\\

\begin{figure}[hbtp]
  \begin{center}
    \leavevmode
    \epsfxsize = 4.5cm
    \epsfysize = 3cm
    \epsfbox{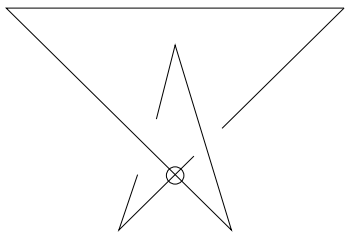}
    \caption{PL representation of $2.1$}
    \label{twoone}
  \end{center}
\end{figure}

\noindent \textit{Case 1}: Suppose $c_{v}(K) = 1$.  Consider the two edges forming the virtual crossing, $e_{1}$ and $e_{3}$, with endpoints the vertices $v_{1}$, $v_{2}$ and $v_{3}, v_{4}$, respectively (all vertices are distinct as $e_{1}$ and $e_{3}$ cannot be adjacent).  These two edges locally separate the plane into four regions and the fifth vertex $v$ of our knot must be placed in one of these regions.  Say it lies in region $R_{1}$ with $R_{i}$ ($i = 2,3,4$) the remaining three regions, counterclockwise from $R_{1}$ around the virtual crossing.  See Fig. \ref{possibles}.

\begin{figure}[hbtp]
  \begin{center}
    \leavevmode
    \epsfxsize = 6.3cm
    \epsfysize = 2.2cm
    \put(25,58){$v$}
    \put(-10,50){$v_{1}$}
    \put(55,-5){$v_{2}$}
    \put(55,50){$v_{3}$}
    \put(-10,-5){$v_{4}$}
    \put(22,40){$R_{1}$}
    \put(0,20){$R_{2}$}
    \epsfbox{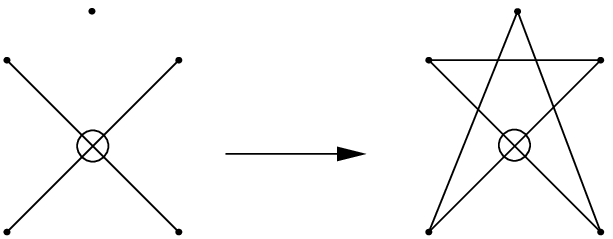}
    \caption{$c_{v}(K) = 1$}
    \label{possibles}
  \end{center}
\end{figure}

We must places the edges $e_{2}$, $e_{4}$, and $e_{5}$ so that we form a knot with crossing index five; that is, so that the diagram is not equivalent to a knot with fewer crossings. Call such a diagram \textit{allowable}. 

Note that one of the edges, say $e_{2}$, must connect an endpoint of $e_{1}$ with an endpoint of $e_{3}$.  If $e_{2}$ is the segment $v_{1}v_{4}$ then $e_{4}$ and $e_{5}$ must necessarily be the segments $v_{3}v$ and $v_{4}v$.  In this case the virtual crossing is removable via a virtual Reidemeister I move, and the knot it unknotted.  Similarly, $e_{2}$ cannot lie in $R_{3}$.   Therefore $e_{2}$ is the segment $v_{1}v_{3}$.  The remaining two edges must be $v_{4}v$ and $v_{2}v$.  This yields a diagram of five crossings.  It is left to determine which choices of real crossings are realizable and yield allowable diagrams.

In order to not have removable crossings, a single edge that passes through two real crossings (in our case, edges $e_{2}$, $e_{4}$, and $e_{5}$) must be the overcrossing strand in one of the crossings and the undercrossing strand in the other, else a Reidemeister II removes these two crossings.  This leaves only two possible diagrams. These diagrams are mirror images of one another and therefore it is enough to consider just one of them.  Notice, however, that this diagram is not realizable, for if it were, let $P$ be the plane in which edges $e_{4}$ and $e_{5}$ lie.  The edge $e_{2}$ crosses through $P$, forcing $v_{1}$ and $v_{2}$ to lie on opposite sides of $P$.  A contradiction arises from how $e_{1}$ and $e_{3}$ lie.  See Fig. \ref{notpossible}.  Thus, there are no knots $K$ with $e(K) = 5$, $c(K) = 5$, and $c_{v}(K) = 1$.\\

\begin{figure}[hbtp]
  \begin{center}
    \leavevmode
    \epsfxsize = 3cm
    \epsfysize = 2.4cm
    \epsfbox{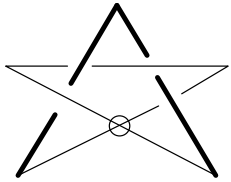}
    \caption{Impossible to form}
    \label{notpossible}
  \end{center}
\end{figure}

\noindent \textit{Case 2}: Suppose $c_{v}(K) = 2$.  The diagram $D$ must have the same form as in the previous case, though we now allow one of the additional four crossings to be virtual.  Via the same edge/crossing stipulations as Case 1 and the fact that a single edge cannot be included in two virtual crossings (else our diagram is not allowable), we are left with the four diagrams in Figure \ref{fivetwovirtual}.

\begin{figure}[hbtp]
  \begin{center}
    \leavevmode
    \epsfxsize = 5.5cm
    \epsfysize = 5.5cm
    \put(27,65){(a)}
    \put(117,65){(b)}
    \put(27,-5){(c)}
    \put(117,-5){(d)}
    \epsfbox{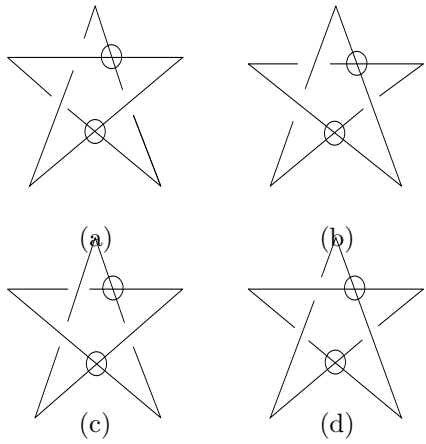}
    \caption{$c_{v}(K) = 2$ possibilities}
    \label{fivetwovirtual}
  \end{center}
\end{figure}

Notice that the pairs (a)-(b) and (c)-(d) in Fig. \ref{fivetwovirtual} are mirror images of one another and we restrict ourselves to two possible knots.  

If $K$ is the knot pictured in Fig. \ref{fivetwovirtual}(a), then $V_{K}(t) = -t^{-4} + t^{-3} + t^{-1}$.  Of the knots in the virtual knot table of crossing index five with two virtual crossings, only $3.5$ has this Jones polynomial.  If $K$ is the knot type of Figure \ref{fivetwovirtual}(c), then $V_{K}(t) = 1$.  There are only two knots of crossing index five and virtual crossing index two in the virtual knot table with this Jones polynomial: $3.1$ and $3.7$.  These two knots are distinguished by the cabled Jones polynomial (see the virtual knot table for its definition).  Consequently we discover that we have formed $3.7$.

Since a PL-virtual knot has the same edge index as its two mirror images, the result follows.
\end{proof}

\subsection{$e(K) = 6$}
\label{equals6}
Consider now six edges.  By Thm. \ref{edgebound}, any virtual knot $K$ realizable with six edges must have $c(K) \leq 8$.  We will restrict our attention to the knots appearing in the virtual knot table (those with $c_{r}(K) \leq 4$).  

\begin{figure}[hbtp]
  \begin{center}
    \leavevmode
    \epsfxsize = 7cm
    \epsfysize = 7.5cm
    \epsfbox{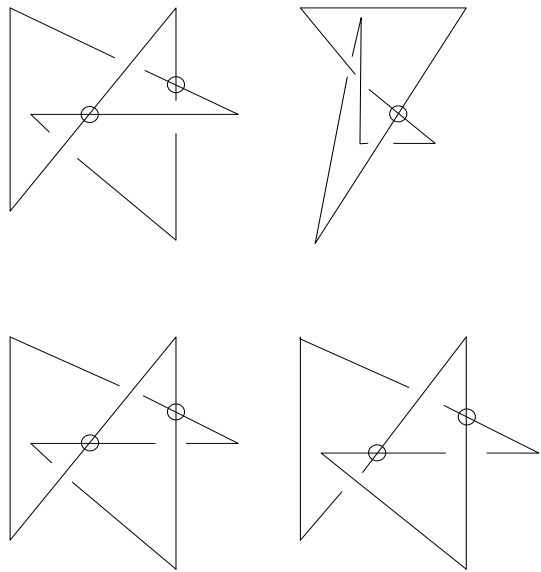}
    \caption{Three crossing virtual knots}
    \label{threecrossings}
  \end{center}
\end{figure}

\begin{lem}
All virtual knots with up to three classical crossings are realizable with six edges.
\end{lem}
\begin{proof}
The knots $0.1$, $2.1$, and $3.5$ were shown to have edge index five in Thm. \ref{fiveedgesorless}.  The remaining four virtual knots are illustrated in Fig. \ref{threecrossings}.
\end{proof}
\vspace{.5cm}

There are over one hundred virtual knots with four real crossings.  Of these, the only possibilities that could be realized by six edges are those with $c_{v}(K) \leq\nolinebreak 3$.  At first glance, it appears that a large number of the virtual knots with $c_{r}(K) = 4$ exceed this virtual crossing bound.  We must be careful, though, as the diagrams in the virtual knot table are only of minimal \textit{real} crossing number and may have equivalent diagrams with fewer virtual crossings.

To proceed, consider the possible shadows of six-edged knots.

Two $n$-universes $U_{1}$ and $U_{2}$ are said to be \textit{isotopic} if there is a continuous deformation $f_{t}$ mapping $U_{1}$ to $U_{2}$ such that at every $t$, $f_{t}(U_{1})$ is in general position (i.e., during the deformation, no vertex crosses an edge).

\begin{lem}
A $6$-universe has at most seven intersections.
\end{lem}
\begin{proof}
This follows from the fact that a $5$-universe has at most five intersections.
\end{proof}

We now consider the isotopy classes of $6$-universes and the knots they arise from.

\begin{figure}[hbtp]
  \begin{center}
    \leavevmode
    \epsfxsize = 7cm
    \epsfysize = 5.5cm
    \put(23,57){($a$)}
    \put(95,57){($b$)}
    \put(165,57){($c$)}
    \put(23,-10){($d$)}
    \put(95,-10){($e$)}
    \put(165,-10){($f$)}
    \epsfbox{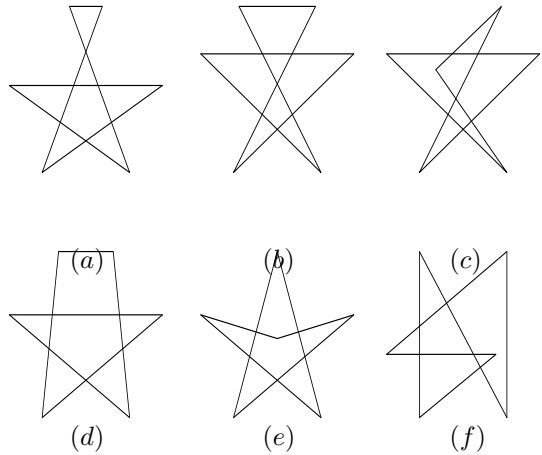}
    \caption{$6$-universes}
    \label{sixuniverses}
  \end{center}
\end{figure}

\begin{lem}
When considering $6$-universes of PL-knots, it is sufficient to consider only one of seven intersections, one of six intersections, and three of five intersections.
\end{lem}
\begin{proof}
Up to isotopy, there is only one such universe of seven intersections.  There are three $6$-universes of six intersections, up to isotopy and mirror images, as pictured in Figs. \ref{sixuniverses}($a$)-($c$).  Note that any knot with the universe in Fig. \ref{sixuniverses}($a$) can be reduced to a diagram with fewer crossings.  Also, any knot realizable with the universe in Fig. \ref{sixuniverses}($c$) is realizable with the universe in Fig. \ref{sixuniverses}($b$).  The three possible $5$-universes are pictured in Figs. \ref{sixuniverses}($d$)-($f$).
\end{proof}

Consider six line segments connected cyclically and lying in the plane.  To be a universe of a virtual knot with four real crossings, the lines must intersect at least five times.  By the previous lemma we know the possible shadows of such diagrams.  Therefore there are only finitely many such PL-virtual knot diagrams that can be formed with six edges.  Enumerate the possible diagrams and using a computer program written by Jeremy Green \cite{VirtualTable} calculate the resulting cabled Jones polynomials and generalized Alexander polynomials for these diagrams.  Similarly, calculate the possible diagrams with five real crossings.  We have the following:

\begin{thm}
The only virtual knots of real crossing index four, up to mirror images, in the virtual knot table realizable with six edges are $4.15$, $4.20$, $4.23$, $4.35$, $4.36$, $4.37$, $4.38$, $4.40$, $4.41$, $4.43$, $4.50$, $4.61$, $4.63$, $4.67$, $4.82$, $4.83$, $4.84$, $4.85$, $4.86$, $4.87$, $4.88$, $4.89$, $4.90$, $4.93$, $4.94$, $4.96$, $4.97$, $4.98$, and $4.99$, as well as all knots realizable with five edges.  There are seventeen distinct knots with real crossing index five that are realizable with six edges.
\end{thm}

\section{Questions}

It is natural to ask if there is a relationship between the edge index of virtual knots and that of classical knots.  In particular, are there bounds that can be placed on the index for virtual knots that in turn strengthen the known bounds for classical knots?  In the spirit of Randell's work on classical PL knots \cite{RandellInvariants}, it would be enlightening to know more about the configuration space of PL virtual knots.  Specifically, is there are natural way to define equivalence in this space?

\bibliography{topology}
\end{document}